\documentstyle[12pt]{article}
\newcommand{\const}{\mathop{\rm const}\limits}

\newcommand{\supp}{\mathop{\rm supp}\limits}

\newcommand{\card}{\mathop{\rm card}\limits}

\newcommand{\Law}{\mathop{\rm Law}\limits}

\newcommand{\grad}{\mathop{\rm grad}\limits}

\newcommand{\argmax}{\mathop{\rm argmax}\limits}

\begin{document}

\begin{center}

{\bf  Relations between exponential tails, moments and \\

\vspace{3mm}

moment generating functions    \\

\vspace{3mm}

for random variables and vectors } \\

\vspace{5mm}

 {\sc Kozachenko Yu.V., \ Ostrovsky E., \ Sirota L.}\\

\vspace{3mm}

 \ Department  of Probability Theory and Statistics,  Kiev State University, \\
Kiev, Ukraine, \\
e-mails:yvk@univ.kiev.ua  \ ykoz@ukr.net \\

\vspace{3mm}

 \ Department of Mathematics and Statistics, Bar-Ilan University, \\
59200, Ramat Gan, Israel.\\
e-mail:eugostrovsky@list.ru \\

\vspace{3mm}

 \ Department of Mathematics and Statistics, Bar-Ilan University,\\
59200, Ramat Gan, Israel.\\
e-mail:sirota3@bezeqint.net \\

\vspace{3mm}

 {\sc Abstract}

\end{center}

 \  We offer in this paper the non-asymptotical pairwise bilateral exact up to multiplicative constants interrelations
between exponential decreasing tail behavior, moments (Grand Lebesgue Spaces) norm and moment generating functions  norm
for random variables and vectors (r.v.). \par

 \vspace{3mm}

 \par

\vspace{4mm}


\vspace{5mm}

\section{Definitions. Statement of problems. Previous results.} \par

\vspace{4mm}

  \hspace{5mm}  Let $ \ (\Omega, F, {\bf P} \ ) $ be a probability space with non - trivial probability measure $ \ {\bf P} $
and expectation $ \ {\bf E}, \ \Omega = \{\omega\}. $  Let also $ \ \xi = \xi(\omega), \ \xi: \Omega \to R $  be numerical valued
random variable (r.v.), i.e. measurable function. The multivariate case may be considered further. \par

\vspace{4mm}

 \ {\bf A. Tail functions.}  \par

 \vspace{4mm}

 \  This function $ T_{\xi}(y), \ y \ge 0 $  for the r.v. $ \ \xi \ $ is defined as ordinary by the formula

$$
T_{\xi}(y) \stackrel{def}{=} {\bf P}(|\xi| \ge y ), \ y \ge 0. \eqno(1.1)
$$

 \ The properties of these functions are obvious. \par

 \ Note that sometimes was used some modifications of these notion, for instance

$$
\overline{T}_{\xi}(y) \stackrel{def}{=} \max[ {\bf P}(\xi \ge y ),  \ {\bf P}(\xi \le - y ) ], \ y \ge 0. \eqno(1.1a)
$$

 \ It is clear

$$
\overline{T}_{\xi}(y) \le T_{\xi}(y) \le 2 \overline{T}_{\xi}(y).
$$

\vspace{4mm}

 \ {\bf B. Grand Lebesgue Spaces.}  \par

 \vspace{4mm}

 \ We recall first of all some needed facts about Grand Lebesgue Spaces (GLS). \par

 \ Recently, see \cite{Kozachenko1}, \cite{Ostrovsky1}, \cite{Ostrovsky5}, \cite{Liflyand1}, \cite{Ostrovsky101},
 \cite{Ostrovsky102}, \cite{Ostrovsky103} etc. appear the so-called Grand Lebesque Spaces $ GLS = G(\psi) =
G(\psi; b), \ b = \const \in (1, \infty] $ spaces consisting on all the random variables (measurable functions)
$ f: \Omega \to R $ with finite norms

$$
     ||f||G(\psi) = G(\psi,b) \stackrel{def}{=} \sup_{p \in [1,b)} \left[ \frac{|f|_p}{\psi(p)} \right]. \eqno(1.2)
$$

  \ Here and in the sequel

$$
|f|_p = [ {\bf E} |f|^p ]^{1/p} = \left[\int_{\Omega} |f(\omega)|^p \ {\bf P}(d \omega)  \right]^{1/p}; \ p \ge 1,
$$
is the classical Lebesgue-Riesz $ L(p) $ norm;
and  $ \psi(\cdot) $ is some continuous positive on the semi - open interval $ [1,b) $ function such that

$$
     \inf_{p \in [1,b)} \psi(p) > 0.
$$

 \ It is proved that $ G(\psi; b) $ is Banach functional rearrangement invariant (r.i.) space and
 $ \supp(G(\psi_b)) := \supp \psi_b = [1,b). $\par

\ Let the {\it family} of measurable functions
$ h_{\alpha} = h_{\alpha}(x), \ x \in X, \ \alpha \in A,  $ where $  \ A \ $ be arbitrary set, be such that

$$
\exists b \in (1, \infty], \ \forall p \in [1, b) \ \Rightarrow
\psi^A(p) := \sup_{\alpha \in A} | \ h_{\alpha} \ |_p < \infty.
$$

 \ Such a function $  \psi^A(p) $ is named as a {\it  natural function} for the family $  A. $ Obviously,

$$
\sup_{\alpha \in A} || \ h_{\alpha} \ ||G\psi^A = 1.
$$

 \ Of course, the family $ \ A  \ $ may consists on the unique function, say, $ \ h = h(\omega), \ $ for which

$$
\exists b > 1 \ \forall p \in [1,b)  \Rightarrow |h|_p < \infty.
$$

   \ These spaces are used, for example, in the theory of probability, theory of PDE, functional analysis,
 theory of Fourier series, theory of martingales etc.\par

\vspace{4mm}

 \ {\bf C. About moment generating function (MGF).}  \par

\vspace{4mm}

 \hspace{3mm} We present here for beginning some known facts from the theory of one-dimensional random variables
with exponential decreasing tails of distributions, see  \cite{Ostrovsky4},  \cite{Kozachenko1}, \cite{Ostrovsky1},
chapters 1,2. \par

 \ Especially  we mention the authors preprint \cite{Ostrovsky5}; we offer in comparison with existing there results
a more fine approach. \par

\vspace{3mm}

 \ Let $ \phi = \phi(\lambda), \lambda \in (-\lambda_0, \lambda_0), \ \lambda_0 =
\const \in (0, \infty] $ be certain even strong convex which takes positive values for positive arguments twice continuous
differentiable function, briefly: Young-Orlicz function, such that
$$
 \phi(0) = \phi'(0) = 0, \ \phi^{''}(0) > 0, \ \lim_{\lambda \to \lambda_0} \phi(\lambda)/\lambda = \infty. \eqno(1.3)
$$

 \ For instance: $  \phi(\lambda) = 0.5 \lambda^2, \ \lambda_0 = \infty; $ the so-called subgaussian case. \par

  \ We denote the set of all these Young-Orlicz  function as $ \Phi; \ \Phi = \{ \phi(\cdot) \}. $ \par

 \ We say by definition that the {\it centered} random variable (r.v) $ \xi = \xi(\omega) $
belongs to the space $ B(\phi), $ if there exists certain non-negative constant
$ \tau \ge 0 $ such that

$$
\forall \lambda \in (-\lambda_0, \lambda_0) \ \Rightarrow
\max_{\pm} {\bf E} \exp(\pm \lambda \xi) \le \exp[ \phi(\lambda \ \tau) ]. \eqno(1.4)
$$

 \ The minimal non-negative value $ \tau $ satisfying (1.4) for all the values $  \lambda \in (-\lambda_0, \lambda_0), $
is named a $ B(\phi) \ $ norm of the variable $ \xi, $ write

$$
||\xi||B(\phi)  \stackrel{def}{=}
$$

 $$
 \inf \{ \tau, \ \tau > 0: \ \forall \lambda:  \ |\lambda| < \lambda_0 \ \Rightarrow
  \max_{\pm}{\bf E}\exp( \pm \lambda \xi) \le \exp(\phi(\lambda \ \tau)) \}. \eqno(1.5)
 $$

 \ These spaces are very convenient for the investigation of the r.v. having a
exponential decreasing tail of distribution, for instance, for investigation of the limit theorem,
the exponential bounds of distribution for sums of random variables,
non-asymptotical properties, problem of continuous and weak compactness of random fields,
study of Central Limit Theorem in the Banach space etc.\par

 \ The space $ B(\phi) $ with respect to the norm $ || \cdot ||B(\phi) $ and
ordinary algebraic operations is a rearrangement invariant Banach space which is isomorphic to the subspace
consisting on all the centered variables of Orlicz's space $ (\Omega,F,{\bf P}), N(\cdot) $
with $ N \ - $ function

$$
N(u) = \exp(\phi^*(u)) - 1, \ \phi^*(u) = \sup_{\lambda} (\lambda u - \phi(\lambda)).
$$
 \ The transform $ \phi \to \phi^* $ is called Young-Fenchel, or Legendre transform. The proof of considered
assertion used the properties of saddle-point method and theorem of Fenchel-Moraux:
$$
\phi^{**} = \phi.
$$

 \ Recall also the Young's inequality

$$
\lambda u \le \phi(\gamma u) + \phi^*(\lambda/\gamma), \ \gamma = \const > 0.
$$

  \ Let $  F =  \{  \xi(s) \}, \ s \in S, \ S  $ is an arbitrary set, be the family of somehow
dependent mean zero random variables. The function $  \phi(\cdot) $ may be "constructive" introduced by the formula

$$
\phi(\lambda) = \phi_F(\lambda) \stackrel{def}{=} \max_{\pm} \ln \sup_{s \in S}
 {\bf E} \exp(  \pm \lambda \xi(s)), \eqno(1.6)
$$
 if obviously the family $  F  $ of the centered r.v. $ \{ \xi(s), \ s \in S \} $ satisfies the  so - called
{\it uniform } Kramer's condition:
$$
\exists \mu \in (0, \infty), \ \sup_{s \in S} T_{\xi(s)}(y) \le \exp(-\mu \ y),\ y \ge 0.
$$
 \ In this case, i.e. in the case the choice the function $ \phi(\cdot) $ by the
formula (1.6), we will call the function $ \phi(\lambda) = \phi_0(\lambda) $
as a {\it natural } function, and correspondingly the function

$$
\lambda \to {\bf E} e^{\lambda \xi}
$$
is named often as a moment generating function for the r.v.  $ \xi, $ if of course there exists  in some non-trivial
neighborhood of origin. \par

 \ Moreover, see \cite{Bagdasarova1},  if $ b = \infty, $ then the following implication holds:

$$
\lim_{\lambda \to \infty} \phi^{-1}(\log {\bf E} \exp(\lambda \xi))/\lambda =
K \in (0, \infty)
$$
if and only if

$$
\lim_{x \to \infty} (\phi^*)^{-1}( |\log T_{\xi}(x)| )/x = 1/K.
$$

 \ If for example $  \phi = \phi_2(\lambda) = 0.5 \ \lambda^2, \ \lambda \in R, $ then the r.v. from the space $ B\phi_2 $ are named
subgaussian. This notion was introduced by  J.P.Kahane  in \cite{Kahane1}; see also \cite{Buldygin2}, \cite{Buldygin3}. One can consider
also the case when

$$
\phi_m(\lambda) = m^{-1} \ |\lambda|^m, \ |\lambda| \ge 1, \ m = \const > 1,
$$
see \cite{Kozachenko1}. \par

\vspace{4mm}

\ {\bf  The aim of this report is to establish the reciprocal non - asymptotic
interrelations separately mutually possibly exact up to multiplicative constant  between tail functions,
 moment generating functions, and Grand Lebesgue Spaces norms.} \par

\vspace{4mm}

 \ {\it Throughout this paper, the letters $ \ C, C_j(\cdot) $ will denote a various
positive finite constants which may differ from one formula to the
next even within a single string of estimates and which does
not depend on the essentially variables $ p, x, \lambda, y $ etc. \par
 We make no attempt to obtain the best values for these constants.} \par

\vspace{4mm}

 \ Obtained here results are unimprovable and generalized ones in  \cite{Kozachenko1}, \cite{Ostrovsky1}, chapter 4;
\cite{Ostrovsky5}, \cite{Ostrovsky6}, \cite{Ostrovsky13}. \par

 \ The applications of these estimates appear for instance in the theory of (discontinuous, in general case) random fields and
following in statistics, see, e.g. in \cite{Bickel1}, \cite{Billingsley1}, \cite{Billingsley2}, \cite{Gikhman1}, chapter 11;
 \cite{Neuhaus1}, \cite{Prokhorov1}, \cite{Prokhorov2}, \cite{Skorokhod1}; in the theory of Monte-Carlo method-in \cite{Frolov1},
\cite{Grigorjeva1}.\par

\vspace{4mm}

\section{Connection between tails and moments}

\vspace{4mm}

 \ {\bf A.  \ "Direct estimate". } Given: the random variable $  \  \xi \ $ such that for some function
 $  \ \psi(\cdot) \in \Psi = \Psi_{\infty} \hspace{4mm} ||\xi|| = ||\xi||G\psi \in (0,\infty). $  It is required to estimate the
 tail function $  \ T_{\xi}(x), $ for sufficiently greatest values $ x, $ say $  \ x > e. \ $\par

 \ Define the auxiliary function

$$
\nu(p) =  \nu_{\psi}(p) := p \ln \psi(p), \ p \ge 1 \eqno(2.1)
$$
and correspondingly

$$
\nu^*(z) = \nu^*_{\psi}(z) = \sup_{p \ge 1} \left( pz - \nu_{\psi}(p) \right), \ z \ge 0. \eqno(2.2)
$$

 \vspace{3mm}

 \ {\bf Theorem 2.1.} Suppose $ \ \xi \in G\psi, \ \xi \ne 0. \ $ Our statement:

$$
T_{\xi}(y) \le \exp \left\{ \ - \nu^*[\ \ln (y/||\xi||) \ ] \ \right\}, \ y > e. \eqno(2.3)
$$

 \vspace{3mm}

\ {\bf Proof.}  \ One can assume without loss of generality $  || \xi ||G\psi = 1. $ Then

$$
|\xi|_p \le \psi(p), \ p \in [1,\infty); \hspace{5mm} \Leftrightarrow \hspace{5mm}
{\bf E} |\xi|^p \le \psi^p(p),  \eqno(2.4)
$$
and we apply the Markov-Tchebychev inequality

$$
T_{\xi}(y) \le \frac{\psi^p(p)}{y^p} =
$$

$$
\exp \left( - (p \ln y - p \ln \psi(p)  \right) =
\exp \left( - (p \ln y - \nu_{\psi}(p)) \right).
$$

 \ It remains to take the minimum over $  p \ge 1 $
 from the right-hand side  of the inequality (2.4); \ $ \ y \ge e. $ \par

\vspace{4mm}

 \ {\bf Remark 2.1.} Let us define  the following Young-Orlicz function

$$
N(u) :=  \exp \left\{ \ \nu^*[ \ \ln |u| \ ] \ \right\}, \ |u| \ge e,
$$
and as usually

$$
N(u): = C \ u^2, \ |u| < e,
$$
where of course  $  C \ e^2 = \nu^*(1). $ \par

 \ It is proved in  \cite{Kozachenko1}, see also \cite{Ostrovsky6}, \cite{Ostrovsky13},
that if the function $ \nu(p) =  p \ \ln \psi(p), \ p \ge 1  $ is continuous, convex,  and such that
$ \lim_{\to \infty} \psi(p) = \infty, $ then the Grand Lebesgue Space $ G\psi $ coincides up to norm equivalence
with Orlicz's space $ L(N) $ builded over source probability space. \par

\vspace{4mm}

 \ {\bf B. "Inverse estimate." } Given:

 $$
 T_{\xi}(x) \le \exp (  - \zeta(x)), \ x \ge 0, \ \zeta(x) \ge 0; \eqno(2.5)
 $$
and we denote $ Z(y) = \zeta(\exp y), \ y \in R.  $ \par

\vspace{3mm}

 \ {\bf Theorem 2.2.}  Suppose $ Z(\cdot) $ is twice continuous differentiable
convex function on certain interval $ \ (C, \infty) \ $  and suppose

$$
\exists C_1 = \const > 0 \ \Rightarrow C_2:= \inf_{y \ge C_1} Z^{''}(y) > 0. \eqno(2.6)
$$
 Then

$$
|\xi|_p \le C_3 \exp \left(  Z^*(p)/p \right), \eqno(2.7)
$$
or equally

$$
||\xi||G(Z^*(p)/p) \le C_3 < \infty. \eqno(2.7a)
$$

\vspace{3mm}

\ {\bf Proof.} We have

$$
{\bf E}|\xi|^p = p \int_0^{\infty} x^{p-1} T_{\xi}(x) dx \le p \ \int_{-\infty}^{\infty} e^{py - Z(y)} dy =: J(p). \eqno(2.8)
$$
 \ It follows from the saddle-point method that

$$
J(p) \le C_4^p \exp \left( \sup_y (py - Z(y)) \right) =  C_4^p \exp \left( Z^*(p) \right).
$$

 \ Let us represent a rigorous consideration. We deduce splitting the integral $  J = J(p) $ onto two ones

$$
J(p) = \ p \ \int_{-\infty}^{ C_1 } e^{py - Z(y)} dy + p \ \int_{C_1}^{\infty}e^{py - Z(y)} dy = J_1(p) + J_2(p).
$$
 \ Note first of all

$$
J_1(p) \le p \ \int_{-\infty}^{ C_1 } e^{py } dy  \le  C_5^p, \ p \ge 1, \ 0 < C_5 < \infty.
$$
 \ Denote $  S = S(p,y) = py - Z(y) $ and $  y_0 = y_0(p) = \argmax_{y \ge C_1} S(p,y);  $ then

$$
S(p,y) \le \max_y S(p,y) - 0.5 S^{''}_y(p,y) (y - y_0)^2 \le Z^*(p) - C_6 (y - y_0)^2,
$$
therefore

$$
J_2(p) \le C_7 \ p \exp \left( Z^*(p) \right)
$$
and following

$$
J(p) \le C_8^p  \exp \left( Z^*(p) \right), \ p \ge 1.
$$

 \ We used the obvious estimate $  p^{1/p} \le C = e^{1/e}. $ \par

 \ This completes the proof of proposition 2.2. \par

\vspace{4mm}

{\bf C. Coincidence.} \par

\vspace{4mm}

 \ It is convenient for us to rewrite the restriction (2.4) in the following form

$$
|\xi|_p \le p \ e^{ \ - \nu(p)/p \ }, \ p \ge 1, \eqno(2.9)
$$
i.e. in (2.4)

$$
\psi(p) := \psi_{\nu}(p) =  \ \ e^{ - \nu(p)/p}.\eqno(2.10)
$$

\vspace{4mm}

 \ {\bf Theorem 2.3.}  Suppose that the function $  \nu = \nu(p), \ p \ge 1  $ in (2.9) is continuous, convex, and
such that the function $  y \to \nu(\exp y) $ satisfies the condition (2.6). Then the GLS norm estimate for the
non - zero r.v. $ \xi $ of the form

$$
|\xi|_p \le C_1 \ p \ e^{ \ - \nu(p)/p \ }, \ p \ge 1, \eqno(2.11)
$$
is quite equivalent to the following tail inequality

$$
T_{\xi}(y) \le \exp \left( - \nu^*(\ln (y/C_2)  \right), \  C_2 = \const \in (0,\infty), \ y \ge C_2 e. \eqno(2.12)
$$

\vspace{4mm}

 \ {\bf Proof.}  The implication (2.11) $ \to $ (2.12) contains really in the statement of theorem 2.1. Conversely, let (2.12) there holds.
It follows from theorem 2.2 that the r.v. $  \xi  $ belongs to the space $  G\psi_{\nu}, \ \xi \in  G(\psi_{\nu}), $ where

$$
\psi_{\nu}(p) = \exp \left(\nu^{**}(p)/p  \right), \ p \ge 1.
$$

 \ But $ \nu^{**} = \nu  $ by virtue of theorem of Fenchel-Moraux, therefore $ \xi \in  G(\psi). $ \par

\vspace{3mm}

 \ This completes the proof of theorem 2.3.\par

\vspace{4mm}

\section{Relations between tails and moments generating functions}

\vspace{4mm}

\ {\bf A. "Direct estimate." } Given as above, see (2.5), for the {\it centered} r.v. $  \ \xi \ $

 $$
 T_{\xi}(x) \le \exp (  - \zeta(x)), \ x \ge 0.
 $$
\ It is required to estimate for the sufficiently greatest values $ \ \lambda, $  say, $ \lambda > e,  $
 the moment generating function (MGF) for the r.v. $  \xi: $

$$
g_{\xi}(\lambda) \stackrel{def}{=} {\bf E} e^{\lambda \xi}, \ \lambda = \const \in R, \eqno(3.1)
$$
or equally the $  \ B(\phi) \ $ norm of the r.v. $  \xi $ for suitable function $ \phi(\cdot) \in \Phi. $ \par
 \ The alternative case $ \lambda \in (-e,e) $ was considered in  \cite{Ostrovsky1}, chapter 1, section 1.2. Note that we
represent here a new approach. \par

\vspace{3mm}

 \ {\bf We must investigate previously  one interest integral.}  Namely, let $ (X,M,\mu) $ be non-trivial measurable space with
non-trivial sigma finite measure $  \mu. $ We assume at once  $ \mu(X) = \infty, $  as long as the opposite case is trivial for us.\par
 \ We intend to estimate for "greatest" values of real parameter $ \lambda, \ \lambda > e $ the following integral

$$
I(\lambda) := \int_X e^{\lambda x - \zeta(x)} \mu(dx) = \int e^{\lambda x - \zeta(x)} \mu(dx), \eqno(3.2)
$$
 assuming of course its convergence for all the sufficiently  great values $ \ \lambda, $ say $ \lambda > e. $\par
Here $ \zeta = \zeta(x) $  is non-negative measurable function, not necessary to be convex. \par

 \ If in contradiction the measure $ \mu $ is finite: $ \mu(X) = M \in (0,\infty); $ then the integral $  I(\lambda) $
allows a simple estimate

$$
I(\lambda) \le M \cdot \sup_{x \in X} \left\{e^{\lambda x - \zeta(x)}\right\} = M \cdot e^{\zeta^*(x)}. \eqno(3.2a)
$$

 \ Let now $ \mu(X) = \infty $  and $ \ \epsilon = \const \in (0,1); $  let us introduce the following integral

$$
K(\epsilon) := \int_X e^{-\epsilon \ \zeta(x)} \ \mu(dx). \eqno(3.3)
$$

 \ It will be presumed its finiteness for all the positive values $ \epsilon > 0; $  or at last for {\it some} positive value
$ \ \epsilon_0 \in (0,1); $  then $ \ \forall \epsilon \ge \epsilon_0 \ \Rightarrow  K(\epsilon) < \infty. $ \par

 \  Then the following measures are probabilistic:

$$
\nu_{\epsilon}(dx) := \frac{e^{-\epsilon \ \zeta(x)} }{K(\epsilon)} \mu(dx): \hspace{6mm} \int_X \nu_{\epsilon}(dx) = 1,
0 < \epsilon < 1. \eqno(3.4)
$$

 \ We have

$$
\frac{I(\lambda)}{K(\epsilon)} = \int \exp(\lambda x - (1 - \epsilon) \ \zeta(x) ) \ \nu_{\epsilon}(dx) \le
$$

$$
\exp \left\{ \ \sup_{x \in X}  [\lambda x - (1 - \epsilon) \ \zeta(x) ] \ \right\} =
\exp \left\{ \ (1 - \epsilon) \sup_x \left[ \ \frac{\lambda}{1 - \epsilon} \ x - \zeta(x) \ \right] \ \right\} =
$$

$$
\exp \left\{ \ (1 - \epsilon) \ \zeta^* \left( \ \frac{\lambda}{1 - \epsilon} \ \right) \ \right\}.
$$
 \ Following,

$$
I(\lambda) \le K(\epsilon) \ \exp \left\{ \ (1 - \epsilon) \ \zeta^* \left( \ \frac{\lambda}{1 - \epsilon} \ \right) \ \right\} \eqno(3.5)
$$
and hence: \par

\vspace{4mm}

{\bf Lemma 3.1} We assert under formulated here conditions:

$$
I(\lambda) \le \inf_{\epsilon \in (0,1)}
 \left[ \ K(\epsilon) \ \exp \left\{ \ (1 - \epsilon) \ \zeta^* \left( \ \frac{\lambda}{1 - \epsilon} \ \right) \ \right\} \ \right].
 \eqno(3.6)
$$

\vspace{4mm}

 \ We intend to simplify the last estimate under some simple additional conditions. In order to carry out this,
 we define the function

$$
\theta(\lambda) \stackrel{def}{=} \frac{C_1}{\lambda \ \zeta^{*'}(\lambda)} \eqno(3.7)
$$
for the greatest values $  \lambda: \ \lambda > \lambda_0, $ where $ \ \lambda_0, \ \lambda_0 = \const, \ \theta(\lambda_0) \le 0.5 $ (say). \par
 \ There is a reasonable to choose in (3.5), (3.6) $ \epsilon := \theta(\lambda), $ see \cite{Ostrovsky1}, chapter 3, pp. 99 - 110.
 \ We conclude denoting

$$
\overline{K}(\lambda) := \int_X e^{ - \theta(\lambda) \ \zeta^*(x) } \mu(dx) = K(\theta(\lambda)):
$$

$$
I(\lambda) \le C_2 \ \overline{K}(\lambda) \ \exp \left( \zeta^*(\lambda)  \right). \eqno(3.8)
$$

\vspace{3mm}

 \ {\bf Definition 3.1.} We will say that the function $  \  \zeta = \zeta(x) \ $ is regular, iff

$$
\exists C_3 = \const < \infty, \ \ \forall \lambda > \lambda_0 \ \Rightarrow \overline{K}(\lambda) \le  \exp \zeta^*(C_3 \lambda). \eqno(3.9)
$$

 \ It follows immediately from lemma 3.1 \\

\vspace{3mm}

{\bf Lemma 3.2.} We assert under formulated here conditions and under the condition of regularity (3.9)

$$
\exists C_4 = \const = C_4(\zeta)  < \infty \ \Rightarrow  I(\lambda) \le \exp \zeta^*(C_4 \lambda), \ |\lambda| > \lambda_0.  \eqno(3.10)
$$

\vspace{3mm}

 \ {\bf We return to the formulated above problem "A". } Indeed, let the estimate (2.5) be a given. Let us estimate the MGF function

$$
g_{\xi}(\lambda) = {\bf E} e^{\lambda \xi}, \ \lambda \ge e.
$$
 \ The case $ |\lambda| \le e $ may  by simple investigated by Taylor's formula. \par

\vspace{3mm}

 \ {\bf Theorem 3.1.} \ Suppose  $  {\bf E}\xi = 0 $ and that in the estimate (2.5) the function $ \zeta = \zeta(x)  $
 satisfies all the conditions of the lemma 3.2, relative the ordinary Lebesgue measure $ \mu(dx) = dx $  and $ X = R, $ in
 particular the condition of regularity. We propose

$$
g_{\xi}(\lambda) = {\bf E} e^{\lambda \xi} \le \exp \left( \ \zeta^*(C \lambda) \ \right), \ \exists  C < \infty, \eqno(3.11)
$$
or equally

$$
|| \ \xi \ ||G\zeta^* = C < \infty. \eqno(3.11a)
$$

\vspace{3mm}

 \ {\bf Proof.} Suppose  $  \xi \ne 0 $ and $ \lambda \ge e; $ the case $ |\lambda| \le e $ may  by simple investigated by the
Taylor's formula taking into account the equality $  \  {\bf E}\xi = 0 \  $
and the case $ \lambda \le - e $ is complete symmetric to the case $  \ \lambda \ge e. $ \par

 \ We have through integration "by parts", see for details \cite{Kozachenko1},

$$
g_{\xi}(\lambda) \le 1 + 2 \lambda \int_{0}^{\infty} \exp(\lambda x - \zeta(x)) \ dx.
$$

 \ We use the statement of the lemma 3.2:

$$
\int_{0}^{\infty} \exp(\lambda x - \zeta(x)) \le \exp \zeta^*( \ C_1 \ \lambda ), \ \lambda \ge e.
$$

 \ The estimate of the form

$$
1 + 2 \lambda \exp \zeta^*( \ C_1 \ \lambda ) \le  \exp \zeta^*( \ C_2 \ \lambda \ ), \ C_2 = \const \in (C_1, \infty), \ |\lambda| \ge e,
$$
is proved in particular in the book \cite{Ostrovsky1}, chapter 1, page 25.\par

\vspace{4mm}

\ {\bf B. "Inverse  estimate." } Given: for the centered r.v. $  \  \xi \ $ the inequality of the form

$$
{\bf E } e^{\lambda \xi} \le e^{\kappa(\lambda)}, \ \lambda \in R, \eqno(3.12)
$$
where  $ \kappa(\cdot) $ is some finite on the whole axis $  R  $  even non-negative function. Let for beginning $  x > 0; $
we apply the famous Chernov's estimate, which follows in turn immediately from the Tchebychev-Markov inequality:

$$
{\bf P}(\xi \ge x) \le \frac{e^{\kappa(\lambda)}}{e^{\lambda x}} = e^{ -(\lambda x - \kappa(\lambda)  ) }, \ \lambda > 0,
$$
 therefore

$$
{\bf P}(\xi > x) \le  e^{ -\sup_{\lambda > 0}(\lambda x - \kappa(\lambda)  ) } = e^{ - \kappa^*(x) }, \eqno(3.13)
$$
 and likewise estimate there holds for the "associate" probability $ {\bf P}(\xi \le - x), \ x > 0. $ Thus,
we deduce under the conditions of the pilcrow {\bf B}\par

\vspace{3mm}

 \ {\bf Theorem 3.2.}

$$
T_{\xi}(x) \le e^{ - \kappa^*(x) }, \ x > 0. \eqno(3.14)
$$

\vspace{4mm}

 \ {\bf Remark 3.1.}  In the article \cite{Gorskikh1} is builded an example of the $ \ \phi(\cdot) \ $ function
from the set $ \Phi $ and the  r.v. $ \ \eta \ $ from the space $ B(\phi)  $ with unit norm:  $  || \ \eta \ ||B(\phi)= 1, $
but for which there exists a deterministic sequence $ \ \{x(n)\}, \ n = 1,2, \ldots  $  tending to infinity
with the following property:

$$
 \exists  C = \const \in (0, \infty), \ \Rightarrow T_{\eta}(x(n)) \ge C \ \exp( - \phi^*(x(n)).
$$
 \ Evidently, the r.v. $ \eta $ is not Gaussian distributed still in the case when $ \phi(\lambda) = 0.5 \ \lambda^2, \ \lambda \in R. $ \par

 \vspace{4mm}

{\bf C. Coincidence.} \par

\vspace{4mm}

 \ {\bf Theorem 3.3.} Suppose that for the certain centered  non - zero r.v. $  \  \xi \ $ there holds the inequality of the form

$$
{\bf E } e^{\lambda \xi} \le e^{\kappa(\lambda)}, \ \lambda \in R, \eqno(3.15)
$$
or equally

$$
|| \ \xi \ ||B(\kappa) \le 1,
$$
where  $ \kappa(\cdot) $ is some finite on the whole axis $  R  $  even non-negative continuous convex
function such that

$$
\kappa(0) = 0, \ \lim_{\lambda \to \infty}  \kappa(\lambda)/\lambda = \infty.\eqno(3.16)
$$

 \ We assert that the equality (3.15) is quite equivalent to the following tail estimate

 $$
\exists K \in (0,\infty) \ \Rightarrow T_{\xi}(x) \le \exp \left( - \kappa^*(x/K)   \right), \ x \ge 0, \eqno(3.17)
 $$
and herewith

 $$
 C_1 K \le ||\xi||B(\kappa) \le C_2 K. \eqno(3.18)
 $$

\vspace{3mm}

\ {\bf Proof.} The implication $ (3.15) \to (3.17)   $ with $  K = 1 $ contains in (3.14). Conversely, let (3.17) be satisfied;
one can take $  K = 1. $ \par
 \ It follows from theorem 3.1 that the random variable. $  \xi  $ belongs to the space $  B(\kappa^{**}): \ \xi \in  B(\kappa^{**}). $
But $ \kappa^{**} = \kappa  $ by virtue of theorem of Fenchel-Moraux, therefore $ \xi \in  B(\kappa). $ \par
 \ The other details may be omitted. \par

\vspace{4mm}

\section{Interrelation between moment generating function and ordinary moments}

\vspace{4mm}

 \ {\bf A. "Direct estimate".} Given: the (centered) r.v.  $  \ \xi \ $ belongs to the certain space $  B(\phi), \ \phi \in \Phi: $

$$
{\bf E} e^{\lambda \xi} \le e^{ \ \phi( \lambda  ||\xi||) \ }, \hspace{5mm} ||\xi|| = ||\xi||B(\phi). \eqno(4.1)
$$

 \ It is required to estimate the GLS norm $  ||\xi||G\psi $ for suitable $ \psi \ - $ function. \par

 \ One can take in (4.1) without loss of generality $ ||\xi|| = ||\xi||B(\phi) = 1, $  so that

$$
{\bf E} e^{\lambda \xi} \le e^{ \ \phi( \lambda) }, \ \lambda \in R.  \eqno(4.1a)
$$

 \ Define the function

$$
\beta_{\phi}(y)= \beta(y) := \phi \left(e^y \right), \ y \in R. \eqno(4.2)
$$

 \ We have using an elementary inequality

$$
z^p \le \left( \frac{p}{e}  \right)^p \cdot e^z, \ z > 0, \ p \ge 1:
$$

$$
{\bf E} |\xi|^p \le  \left( \frac{p}{e \lambda}  \right)^p e^{\phi(\lambda)}, \ \lambda > 0,
$$
therefore

$$
  {\bf E}|\xi|^p \le
  \left( \frac{p}{e}  \right)^p \cdot \inf_{\lambda > 1} \left[ e^{ - p \ln \lambda + \phi(\lambda)} \right] =
\left( \frac{p}{e}  \right)^p \cdot \exp \left( - \sup_{\lambda > 1}(p \ln \lambda - \phi(\lambda) \right) =
$$

$$
\left( \frac{p}{e}  \right)^p \cdot  \exp \left( - \sup_{\mu > 0} (p \mu - \beta(\mu)  \right) =
\left( \frac{p}{e }  \right)^p \cdot \exp (- \beta^*(p)),
$$

 \vspace{3mm}

 Result: \par

\vspace{3mm}

 \ {\bf Theorem 4.1.} Denote

$$
\psi_{(\phi)}(p) = p \cdot \exp \left[ - \beta^*(p)/p  \right], \ p \ge 1; \eqno(4.3)
$$
then we have under formulated before conditions and notations

$$
|| \ \xi \ ||G\psi_{(\phi)} \le e^{-1} \ ||\ \xi \ ||B(\phi). \eqno(4.4)
$$

\vspace{4mm}

 \ {\bf B. "Inverse estimate".} Given: the non-zero (centered) r.v.  $  \ \xi \ $ belongs to the certain space
$  G(\psi), \ \psi \in \Psi:  \ ||\xi||G\psi \in (0,\infty). $ We need to estimate the MGF function for the r.v. $  \xi. $ \par

 \ We can and will suppose $ ||\xi||G\psi = 1,  $ therefore

 $$
 |\xi|_p \le \psi(p), \ p  \ge 1.  \eqno(4.5)
 $$

 \ It is convenient for us to rewrite the restriction (4.5) in the following form

$$
|\xi|_p \le p \ e^{ - \Delta(p)/p}, \ p \ge 1, \eqno(4.5a)
$$
i.e. in (4.5)

$$
\psi(p) := \psi_{\Delta}(p) = p \ e^{ - \Delta(p)/p}.\eqno(4.5b)
$$
 \ The concrete conditions on the function $ \Delta(\cdot) $ will be clarified below.\par

 \ Note that the condition

$$
\lim_{p \to \infty} [\psi(p)/p]  = 0
$$
is necessary for the existence of all the exponential moments of the r.v. $ \xi, $ on the other hands, for the existence
of MGF for r.v. $  \ \xi. $ \par

 \ It is enough to use the inequality (4.5) only for the even number $ p: p = 2 m, \ m = 0,1,2, \ldots:  $

$$
{\bf E} \xi^{2m} \le \psi^{2m}(2m) = (2m)^{2m} \exp( - \Delta(2m) ).
$$

 \ We have for the great values $ \lambda, $ say for $ \ \lambda \ge e,  $ using the Stirling's formula

$$
{\bf E}\cosh (\lambda \xi) - 1 =  \sum_{m=1}^{\infty} \frac{\lambda^{2m}}{(2m)!} {\bf E}\xi^{2m} \le
$$

$$
\sum_{k = 2,4,\ldots} \exp( k \ln(C\lambda) - \Delta(k)). \eqno(4.6)
$$

  \ Let $  \epsilon = \const \in (0,1); $ we apply the Young's inequality, denoting $ \mu = \ln (C \lambda): $

$$
k \mu \le \Delta(\epsilon k) + \Delta^*(\mu/\epsilon)
$$
and denoting

$$
\sigma_{\Delta}(\epsilon) = \sigma(\epsilon) := \sum_{k = 2,4,\ldots} e^{ \Delta(\epsilon k) - \Delta(k)}: \eqno(4.7)
$$

$$
{\bf E}\cosh (\lambda \xi) - 1 \le e^{\Delta^*(\mu/\epsilon)} \sum_{k = 2,4,\ldots} e^{ \Delta(\epsilon k) - \Delta(k)}  \le
$$

$$
\sigma_{\Delta}(\epsilon) \cdot e^{\Delta^*(\mu/\epsilon)}.\eqno(4.8)
$$

 \ Introduce also the function $  \Delta_1^*(\mu) $ as follows:

$$
 \Delta_1^*(\mu) := \inf_{\epsilon \in (0,1)} \ln \left\{ \sigma(\epsilon) \ \exp \left( \Delta^*(\mu/\epsilon)  \right) \right\};
 \eqno(4.9)
$$
then we obtain for the values $ \lambda \ge e: $

$$
{\bf E} e^{\lambda \xi} \le e^{ \Delta_1^*(\ln (C_2 \lambda) )}. \eqno(4.10)
$$

 \ At the same estimate is true also with another but again finite constant $ C_3 $ for the values  $ |\lambda| < e $
and $ \lambda \le - e, $ following

$$
 \exists C_3 < \infty \ \Rightarrow  {\bf E} e^{\lambda \xi} \le e^{ \Delta_1^*(\ln (C_3 \ | \lambda)| )}. \eqno(4.11)
$$

\vspace{3mm}

 \ Let us impose the following important condition $ (\Delta) $ on the source function $  \phi_{\Delta}(\cdot):  $

$$
 (\Delta): \hspace{5mm}   \exists C_4 \in [1, \infty) \Rightarrow \Delta_1^*(\ln (\ | \lambda)| ) \le
 \Delta^*(\ln (\ | C_4 \lambda)| ), \ \lambda \in R. \eqno(4.12)
$$

\vspace{4mm}

 \ {\bf Theorem 4.2.} Suppose the function $ \psi = \psi_{\Delta}  $  in (4.5a), (4.5b) satisfies the
condition $ \ (\Delta) \ $ (4.12). Define the following function

$$
\phi_{\Delta}(\lambda) = \Delta^*(\ln |\lambda|), \ \lambda \in R. \eqno(4.13)
$$

 \ Our proposition:

$$
|| \ \xi \ ||B \left(\phi_{\Delta} \right) \le C_4 \ || \ \xi \ ||G\psi_{\Delta}. \eqno(4.14)
$$

\vspace{4mm}

\ {\bf C. Coincidence.}

\vspace{4mm}

 \ {\bf Theorem 4.3.} Let the initial function $ \phi $ from the set $  \Phi $ be such that
the correspondent function $  \Delta_{\phi}:= \beta^*_{\phi} $ satisfies the condition $ (\Delta). $ Then both the norms
$ || \ \xi \ ||B(\phi) $ and $  \ || \ \xi \ ||G\psi_{\phi} $  are equivalent:

$$
|| \ \xi \ ||B (\phi) \le C_4 \ || \ \xi \ ||G\psi_{\phi} \le
C_5 || \ \xi \ ||B (\phi). \eqno(4.15)
$$

\vspace{4mm}

\section{The case of bounded support}

\vspace{4mm}

 \ We investigate in this section the case when the generating functions $  \psi(\cdot), \ \phi(\cdot) $ have a bounded support:
 $ \supp \psi(p) = [1,b) $ or equally $ \supp \psi(p) = [1,b], \  b = \const \in (1,\infty). $
The case $ \supp \phi (\cdot) = (-1,1) $ will be considered othe second half. \par

\vspace{4mm}

 \ {\bf A. \ Tail-moment relations.} \par

\vspace{3mm}

 \ Suppose the r.v. $ \ \xi \ $ is such that

 $$
  \ T_{\xi}(x) \ \le  T^{\beta,\gamma,L}(x/K), \ \beta = \const > 1, \ \gamma = \const  > -1, \ K = \const > 0, \eqno(5.1)
 $$
$  L = L(x) $ is positive slowly varying at infinity continuous function,

$$
T^{\beta,\gamma,L}(x) \stackrel{def}{=} x^{-\beta} \ (\ln x)^{\gamma} \ L(\ln x), \ x \ge 1. \eqno(5.2)
$$
 \ We have for the values $ p:  \ 1 \le p < \beta, \ p \to \beta - 0  $ and taking for the sake of simplicity  $  \ K = 1\ $

$$
{\bf E}|\xi|^p \le C p \int_1^{\infty}x^{p - \beta - 1} \ (\ln x)^{\gamma} \ L(\ln x) \ dx  =
$$

$$
C \ p \ \int_0^{\infty} e^{ -y(\beta - p) } \ y^{\gamma} \ L(y) \ dy = C \ p \ \ (\beta - p)^{-\gamma - 1} \
\int_0^{\infty} e^{-z} \ z^{\gamma} \ L \left(  \frac{z}{\beta - p}  \right) \ dz \sim
$$

$$
C \ p \ \ (\beta - p)^{-\gamma - 1} \  L \left(  \frac{1}{\beta - p}  \right) \
\int_0^{\infty} e^{-z} \ z^{\gamma} \ dz  =
$$

$$
C \ p \ (\beta - p)^{-\gamma - 1}
 \ L(1/(\beta - p)) \ \Gamma(\gamma + 1),
$$
therefore

$$
|\xi|_p \le C_1(\beta,\gamma, L)  \ (\beta - p)^{-(\gamma +1)/\beta} \   L^{1/\beta}(1/(\beta -  p)), 1 \le p < \beta. \eqno(5.3)
$$

 \ Conversely, assume that for certain r.v. $ \ \xi \ $ the estimate (5.3) there  holds. We find as before, see theorem 2.1:

$$
T_{\xi}(x) \le C_2(\beta,\gamma, L) \ x^{-\beta} \ (\ln x)^{\gamma+1} \ L(\ln x), \ x \ge 1. \eqno(5.4)
$$

\vspace{3mm}

 \ {\bf Remark 5.1.} Notice that there is a gap between estimates (5.2) and (5.4). On the other words, in the considered case
 $  b < \infty $  there is not quite coincidence between tail behavior and Grand Lebesgue Space norm estimates. \par

\vspace{3mm}

 \ {\bf Remark 5.2.}  Obviously, in the considered above example the moment generating function does not exists. \par

\vspace{3mm}

\ {\bf Remark 5.3.}  Both the inequalities (5.3) and (5.4) are non-improvable, see
\cite{Liflyand1}, \cite{Ostrovsky102}. \par

\vspace{4mm}

 \ {\bf B. \ Tail-moment generating function relations.} \par

\vspace{3mm}

 \ Let $ \theta = \const > -1,  \  L = L(x) $ be again positive slowly varying at infinity continuous function.\par

 \vspace{4mm}

 \ {\bf Theorem 5.1.} \\

 \vspace{3mm}

{\bf I.} The inequality of the form

$$
T_{\xi}(x) \le C_1 \ x^{\theta} \ L(x) \ e^{-x}, \ x \ge 1 \eqno(5.5)
$$
follows the following  MGF estimate

$$
\phi_{\xi}(\lambda) \le C_2 \ (1 - |\lambda|)^{-1 - \theta} \ L \left( \frac{1}{1 - |\lambda|} \right), \
|\lambda| < 1. \eqno(5.6)
$$
$ C_2 = C_2(C_1,\theta,L).$ \\

 \vspace{3mm}

{\bf II.} Conversely, it follows from the estimate (5.6)

$$
T_{\xi}(x) \le C_3(C_2, \theta,L)  \ x^{\theta + 1} \ L(x) \ e^{-x}, \ x \ge 1. \eqno(5.7)
$$

 \vspace{3mm}

{\bf Proof.} Let the inequality (5.5) be a given. Let also $ \ \lambda \in (0,1), \ \lambda \to 1 - 0. $ The opposite case
$ \lambda \in (-1,0) $ may be investigated analogously. \par
 \ We deduce

$$
C_4^{-1} g_{\xi}(\lambda) \le \int_0^{\infty} e^{ -x(1 - \lambda) } \ x^{\theta} \ L(x) \ dx =
(1 - \lambda)^{-\theta - 1} \int_0^{\infty} e^{-y} \ y^{\theta} \ L \left( \frac{y}{1 - \lambda} \right) \ dy \sim
$$

$$
(1 - \lambda)^{-\theta - 1} \  L \left( \frac{1}{1 - \lambda} \ \right)
\int_0^{\infty} e^{-y} \ y^{\theta} \ dy  =
\Gamma(\theta + 1) \ (1 - \lambda)^{-\theta - 1} \  L \left( \frac{1}{1 - \lambda} \right).
$$

 \vspace{3mm}

 \ Conversely, let (5.6) be given. We apply the proposition of theorem 3.2. In detail, for the positive values $  x  $

$$
\phi_{\xi}^*(x) =
\sup_{\lambda \in (0,1)} \left[ \lambda x + (\theta + 1) \ln(1 - \lambda) + \ln L \left(  \frac{1}{1 - \lambda}  \right) \right].\eqno(5.8)
$$

  \ We get to the required estimate (5.7) by selection in (5.8) the value

$$
\lambda := \lambda_0 \stackrel{def}{=} 1 - \frac{\theta + 1}{x}, \ x \ge 2(\theta +1).
$$

\vspace{4mm}

 \ {\bf C. \ Tail-Grand Lebesgue Spaces norm relations.} \par

\vspace{3mm}

 \ Let the r.v. $  \xi $ satisfies the tail estimate (5.5); then we have for sufficiently greatest valued $  x  $

$$
T_{\xi}(x) \le C_5 e^{-x/2}, \ x \ge C_6. \eqno(5.9)
$$

 \ It follows immediately from (5.9)

$$
\sup_{p \ge 1} \left[ \frac{|\xi|_p}{p}  \right] < \infty, \eqno(5.10)
$$
and conversely it follows from (5.10) the tail estimate

$$
T_{\xi}(x) \le C_7 e^{-x/C_8}, \ x \ge C_9. \eqno(5.11)
$$

\vspace{4mm}

{\bf Remark 5.4.} Introduce the following $ \psi \ -  $ function

$$
\psi_1(p) = p, \ p \in [1, \infty).
$$
 \ The relation (5.10) may be rewritten on the language GLS spaces as follow assertion

$$
||\xi||G\psi_1 < \infty,
$$
witch is in turn equivalent to the relation (5.5). \par

 \ On the other words, in this case both the estimates are equivalent in the sense of belonging of the r.v. $ \xi $
 to at the same Banach space $ \ G\psi_1, \ $ in contradiction to the estimates (5.2)-(5.4). \par

\vspace{4mm}

\section{Some examples and counterexamples}

\vspace{4mm}

\ {\bf Example 1.}

\vspace{4mm}

 \ Let $  L = L(\lambda)  $ be positive twice continuous differentiable slowly varying at infinity regular
in the following sense

$$
\lim_{\lambda \to \infty} \left\{ \frac{L(\lambda/L(\lambda))}{L(\lambda) } \ \right\} = 1
$$
function. Define also for sufficiently greatest values $ \lambda, $ say for $ |\lambda| \ge e, $  the $  \ \Phi \ - \ $ function
$ \phi_{m,L} = \phi_{m,L}(\lambda), \ m = \const \in(1, \infty) $ of the form

$$
\phi_{m,L}(\lambda) \stackrel{def}{=} m^{-1} \ |\lambda|^m  \ L^{1/q}(|\lambda|^m), \ q = m/(m-1),
$$
and as usually

$$
\phi_{m,L}(\lambda) \stackrel{def}{=} C_{m,L} \ \lambda^2, \ |\lambda| \le 1. \eqno(6.1)
$$
 \ The correspondent $ B(\phi) $ space will be denoted by $ B_{m,L} := B(\phi_{m,L}). $ \par

 \ Define  also the following $ \psi \ - $ function

$$
\psi_{m,L}(p) = p^{1/m} L^{-1/(m-1)} \left( \ p^{ (m-1)^2/m } \  \right), \ p \ge 1. \eqno(6.2)
$$

\vspace{3mm}

 \ Let $  \ \xi \ $ be non-zero centered: $ {\bf E} \xi = 0 $ r.v.
 We conclude by virtue of theorem 3.3 that the including  $  \xi \in B_{m,L},  $
 i.e. the inequality

$$
\exists C_1 \in (0,\infty), \ \forall \lambda \in R \ \Rightarrow
 {\bf E}\exp{\lambda \xi} \le
\exp\left( \phi_{m,L}(C_1 \ \lambda)   \right), \ 0 < C_1 < \infty,  \eqno(6.3a)
$$
or equally

$$
|| \ \xi \ ||B\phi_{m,L} \le C_1 < \infty,
$$

is quite equivalent to the following tail estimate

$$
\exists C_2 \in (0,1), \ \forall y \ge e \ \Rightarrow
$$

$$
T_{\xi}(y) \le \exp \left\{ -C_2 \ q^{-1} \ y^q \ L^{ - (q-1)} \ \left(y^{q-1} \ \right) \ \right\}, \ q = m/(m-1), \eqno(6.3b)
$$

or in turn is equivalent to the Grand Lebesgue Space norm estimate

$$
|| \ \xi  \ || G\psi_{m,L} = C < \infty. \eqno(6.3c)
$$

\vspace{3mm}

 \ It is sufficient to verify this statement to mention the book \cite{Seneta1}, pp. 32-33, where is in particular
calculated the Young-Fenchel transform for the  function $ \phi_{m,L}(\lambda). $ \par

 \vspace{3mm}

 \ Let us consider the case when in addition

$$
L(\lambda) = L_r(\lambda) =  [\ln  \ | \lambda|]^r, \ \lambda \ge e, \ r = \const \in R.
$$

\ The correspondent $ \phi \ - $ function will be denoted by $ \phi_{m,r}(\lambda): $

$$
\phi_{m,r}(\lambda)= |\lambda|^m \ (\ln |\lambda|)^r, \ |\lambda| \ge e,
$$

and the correspondent $ \psi(p)  = \psi_{m,r}(p) \ - $ function  has a form

$$
\psi_{m,r}(p)  = p^{1/m} \ \ln^{-r/(m-1)} (p+1), \ p \ge 1.  \eqno(6.4)
$$

\ Introduce also the following tail function

$$
T^{(m,r)}(x) = \exp \left\{ \ - x^q \ (\ln x)^{-(q - 1) r} \ \right\}, \ x \ge e. \eqno(6.5)
$$

 \ The following propositions for the non-zero centered r.v. $  \ \xi \ $ are equivalent:

$$
{\bf A.} \hspace{5mm} || \ \xi \ ||G\psi_{m,r} < \infty; \eqno(6.6)
$$

$$
{\bf B.} \hspace{5mm} ||\xi ||B(\phi_{m,r}) < \infty; \eqno(6.7)
$$

$$
{\bf C.} \hspace{5mm} \exists K = \const \in (0,\infty) \Rightarrow
T_{\xi}(x) \le T^{(m,r)}(x/K); \eqno(6.8)
$$

and herewith both the norms $ \ || \ \xi \ ||G\psi_{m,r},  \ |\xi ||B(\phi_{m,r}) \  $ and the
"Tail constant" $  K  $ are equivalent:

$$
 K \le C_1 \ || \ \xi \ ||G\psi_{m,r} \le C_2 || \ \xi \ ||B(\phi_{m,r})  \le C_3 K, \eqno(6.9)
$$
if we understood as the capacity of the value $  \ K \ $ its maximal value. \par
 \ The last statement generalized ones obtained in \cite{Kozachenko1}, \cite{Neuhaus1},
 \cite{Ostrovsky1},  section 1, pp. 22-26; \cite{Ostrovsky6}, \cite{Ostrovsky204}, \cite{Ostrovsky205} etc.\par

\vspace{4mm}

\ {\bf Example 2.}

\vspace{4mm}

 \ Another example: define the other Grand Lebesgue Space space of random variables $ \Psi(C,\beta), \
\beta = \const > 0 $ which consist on all the random variables $ \{\eta\} $ with finite norm ( $ C \in
(0,\infty) $ )

$$
|||\eta|||_{C,\beta} \stackrel{def}{=} \sup_{p \ge 1} \left[ \ |\eta|_p \exp \left(-C \ p^{\beta} \right)\ \right].
$$
 \ It is easy to verify that $ \eta  \in \cup_{C > 0} \Psi(C,\beta), \eta \ne 0  \Leftrightarrow $

$$
\ T_{\eta}(x) \le \exp \left( - C_1(C,\beta) (\log (1+x))^{1+1/\beta} \right), \ x \ge 0.
$$

 \ Note that in this case the MGF for arbitrary r.v. with

$$
\ T_{\eta,x} \ge \exp \left( - C_1(C,\beta) (\log (1+x))^{1+1/\beta} \right).
$$
does not exists; on the other words this variable does not satisfy the Kramer's condition. \par

 \ Let us represent  more exact conclusions. Define as ordinary for every constant $ \theta > 1  $

$$
\theta' \stackrel{def}{=} \frac{\theta}{\theta - 1}.
$$
 \ Obviously, $  (\theta')' = \theta. $\par

 \ Suppose for certain r.v. $  \eta $

$$
\ln T_{\eta}(y) \le C_1 - \frac{\ln^{\theta} y}{\theta}, \ y \ge e,\eqno(6.10)
$$
then

$$
\ln |\eta|_p \le C_2(C_1,\beta)  +  \frac{p^{\theta'}}{\theta'}, \ p \ge 1,
\eqno(6.10a)
$$
and conversely, the relation (6.9) follows in turn from (6.10). \par

 \ Moreover, the following "Tauberian" conclusion holds true, see  \cite{Bagdasarova1}. The following propositions
are completely equivalent:

$$
\lim_{y \to \infty} \left[  \left\{ \ |\ln T_{\eta}(y)| \ \right\} : \left\{ \frac{\ln^{\theta} y}{\theta} \right\} \right] = 1, \eqno(6.11)
$$

$$
\lim_{p \to \infty}  \left[  \left\{ \ \ln |\eta|_p \ \right\} : \left\{ \ \frac{p^{\theta'}}{\theta'} \ \right\}  \right]  = 1.
\eqno(6.12)
$$

\vspace{4mm}

\ {\bf Counterexample 1.}

\vspace{4mm}

 \ Let $ \ (\Omega, F, {\bf P} \ ) $ be the classical probability space $ ((0,1), F, \mu), $  where
$ \mu(d \omega) = d \omega $ is ordinary Lebesgue measure. Let also $ \ \alpha = \const \in (0, 1/2); \ p_0 := 1/\alpha > 2.  $  \par
 \ We consider here the random variable

$$
\xi(\omega) := \omega^{-\alpha}.
$$

 \ We have:

$$
\psi_{\xi}(p) =  |\xi|_p  = (1 - \alpha p)^{-1/p} \asymp (p_0 - p)^{-\alpha}, \ p \in [1, p_0),
$$

$$
\psi_{\xi}(p) =  + \infty, \ p \ge p_0;
$$

$$
T_{\xi}(y) = y^{-1/\alpha}, \ y \in (0,1);
$$

$$
\nu_{\xi}(p) \asymp |\ln (p_0 - p)|, \ p \in [1, p_0);
$$

$$
\nu^*( \ln y) \asymp -1  + p_0 \ln y - \ln \ln y, \ y > e^e.
$$

 \ Theorem 2.1, more precisely,  the inequality (2.3) gives  us the following estimate

$$
T_{\xi}(y) \le C(\alpha) \ y^{-1/\alpha} \ \ln y, \ y \ge e,
$$
which is "essentially" greatest as $ \ y \to \infty  $ than the exact value of this tail function $ T_{\xi}(y). $ \par

  \ The reason for this phenomenon is that here $ \psi_{\xi}(\cdot) \in G\psi_b, $ where the value $ b $ is finite:
 $  b = p_0 = 1/\alpha < \infty, $  in contradiction to the conditions of theorem 2.1. \par
  \ Many other interest examples about this relations may be found in \cite{Ostrovsky102}.

\vspace{4mm}

\ {\bf Counterexample 2.}

\vspace{4mm}

 \ Let $ \xi $ be symmetrical distributed r.v. having the Laplace distribution with the density

$$
f_{\xi}(x) = f(x) = 0.5 \ e^{-|x|}, \ x \in R.
$$
 \ It is easily to calculate:

$$
\psi_{\xi}(p) = |\xi|_p \asymp p, \ 1 \le p < \infty,
$$

$$
{\bf E} e^{\lambda \ \xi} = \frac{1}{1 - \lambda^2}, \ |\lambda| < 1,
$$
so that

$$
\phi_{\xi}(\lambda) = - \ln (1 - \lambda^2),  \ \ |\lambda| < 1,
$$
and

$$
\phi_{\xi}(\lambda) = +\infty
$$
otherwise;

$$
T_{\xi}(y) = e^{-y}, \ y \ge 0.
$$

 \ The applying of theorem 2.1 gives the estimate

$$
T_{\xi}(y) = e^{-y/C}, \ y \ge Ce, \ C = \const > 1;
$$
 i.e. despite that  in this case again $ \ b = 1 < \infty, $ the proposition of theorem 2.1 remains true. \par

\vspace{4mm}

\section{Multivariate case}

\vspace{4mm}

 \hspace{3mm} The theory of the multidimensional $  B(\phi), \ G\psi $ spaces and the spaces of random vectors with exponential
 decreasing tail if distribution is described in the recent article \cite{Ostrovsky6}; it is quite analogous to the explained one.
 We represent further briefly in this section some generalizations of previous results into the case when the instead the random variable
 stands a {\it random vector,} which is abbreviated also by r.v. \par

 \ "Briefly"-in the case when the multidimensional version is completely alike to the one-dimensional one. \par

 \vspace{3mm}

 \ In detail,  denote by $  \epsilon  = \vec{\epsilon} = ( \epsilon(1),   \epsilon(2), \ldots,   \epsilon(d) ) $ the non-random
$  d \ - $ dimensional  numerical vector, $  d = 2,3,\ldots, $ whose components take the values $  \pm 1 $ only.
 Set in particular $  \vec{1} = (1,1,\ldots,1) \in R^d_+. $ \par

 \ Denote also by $  \Theta = \Theta(d) = \{ \ \vec{\epsilon} \ \} $ {\it collection} of all such a vectors. Note that
 $  \card \Theta = 2^d $  and $ \vec{1} \in \Theta.$  \par

 \ Another notations. For $ \vec{\epsilon} \in \Theta(d) $  and vector $ \vec{x} $ we introduce the coordinatewise tensor product
as a $ d  \ - $ dimensional vector of the form

$$
\vec{\epsilon} \otimes \vec{x} \stackrel{def}{=} ( \epsilon(1) \ x(1), \ \epsilon(2) \ x(2), \ \ldots,  \epsilon(d) \ x(d) ),
$$
and analogously may be defined recursively the triple tensor product

$$
\vec{\epsilon} \otimes \vec{x} \otimes \vec{y} = (\vec{\epsilon} \otimes \vec{x}) \otimes \vec{y} =
$$

$$
  ( \epsilon(1) \ x(1) \ y(1), \ \epsilon(2) \ x(2) \ y(2), \ \ldots,  \epsilon(d) \ x(d) \ y(d) ).
$$

\vspace{3mm}

 \ {\bf Definition 7.1.}\par

\vspace{3mm}

 \ Let $ \xi = \vec{\xi} = (\xi(1), \xi(2), \ldots, \xi(d) )  $ be a centered random vector such that each its component
$  \xi(j) $  satisfies the Kramer's condition. The {\it natural function} $ \phi_{\xi}= \phi_{\xi}(\lambda), \
 \lambda = \vec{\lambda} = (\lambda(1), \lambda(2), \ldots, \lambda(d))  \in R^d  $
for the random vector $  \xi  $ is defined as follows:

$$
\exp \{\phi_{\xi}(\lambda)\} \stackrel{def}{=} \max_{\vec{\epsilon} \in \Theta} \
{\bf E} \exp \left\{ \sum_{j=1}^d  \epsilon(j) \lambda(j) \xi(j) \right\}=
$$

$$
\max_{\vec{\epsilon} \in \Theta} \
{\bf E} \exp \{\epsilon(1) \lambda(1) \xi(1) + \epsilon(2) \lambda(2)\xi(2)+ \ldots +   \epsilon(d) \lambda(d)\xi(d) \}=
$$

$$
\max_{\vec{\epsilon} \in \Theta} \ {\bf E} \exp ( \vec{\epsilon} \otimes \vec{\lambda} \otimes \vec{\xi}).
 \eqno(7.1)
$$
 where $  "\max" $ is calculated over all the combinations of signs $ \epsilon(j) =  \pm 1. $\par

\vspace{3mm}

 \ {\bf Definition 7.2.}\par

\vspace{3mm}

 \ The {\it tail function } for the random vector $ \vec{\xi} \ $
   $ U(\vec{\xi}, \vec{x}), \ \vec{x} = (x(1), x(2), \ldots, x(d)),   $ where all the coordinates $ x(j)  $
  of the deterministic vector $ \vec{x} $ are non-negative, is defined as follows.

$$
U(\vec{\xi}, \vec{x}) \stackrel{def}{=} \max_{\vec{\epsilon} \in \Theta}
{\bf P} \left( \cap_{j=1}^d  \{  \epsilon(j) \xi(j) > x(j) \}  \right) =
$$

$$
 \max_{ \vec{\epsilon} \in \Theta}
{\bf P}(\epsilon(1) \xi(1) > x(1), \ \epsilon(2) \xi(2) > x(2), \ \ldots, \ \epsilon(d) \xi(d) > x(d) ), \eqno(7.2)
$$
 where as before $  "\max" $ is calculated over all the combinations of signs $ \epsilon(j) =  \pm 1. $\par

 \ We illustrate this notion in the case $  d = 2.  $ Let $ \vec{\xi} = (\xi(1), \ \xi(2))  $ be a two-dimensional random
vector and let $ x, y  $  be non-negative numbers. Then for all the non-negative values $ \ x,y  $

$$
U( (\xi(1), \xi(2)), \ (x,y) ) =
$$

$$
\max [ {\bf P} (\xi(1) > x, \ \xi(2) > y), \ {\bf P} (\xi(1) > x, \ \xi(2) < - y),
$$

$$
{\bf P} (\xi(1) < - x, \ \xi(2) > y), \ {\bf P} (\xi(1) < - x, \ \xi(2) < - y) ].
$$

\vspace{3mm}

 \ {\bf Definition 7.3.}\par

\vspace{3mm}

 \ Let $  h = h(x), \ x \in R^d  $ be some non-negative real valued function, which is finite on some non-empty
neighborhood of origin. We denote as ordinary

$$
\supp h = \{x, \ h(x) < \infty   \}.
$$

 \ The Young-Fenchel, or Legendre transform $  h^*(y), \ y \in R^d  $ is defined likewise the one-dimensional case

$$
h^*(y) \stackrel{def}{=} \sup_{x \in \supp h} ( (x,y) - h(x)). \eqno(7.3)
$$

 \ Herewith $ (x,y)  $ denotes the scalar product of the vectors $ x,y: \ (x,y) = \sum_j x(j)y(j);  \ |x| = \sqrt{(x,x)}. $\par

 \ Obviously, if the set $ \supp h $ is central symmetric, then the function $ h^*(y) $ is even.\par

\vspace{3mm}

 \ {\bf Definition 7.4.}\par

\vspace{3mm}

 \ Recall, see \cite{Krasnosel'skii1}, \cite{Rao1}, \cite{Rao2} that the function $  x \to g(x), \ x \in R^d, \ g(x) \in R^1_+  $ is named
 multivariate Young, or Young-Orlicz function, if it is even,  $ d \ - $ times  continuous differentiable, convex, non - negative,
 finite on the whole space $  R^d,$ and such that

$$
g(x) = 0 \ \Leftrightarrow x = 0; \hspace{4mm} \frac{\partial g}{\partial x}/(\vec{x} = 0) = 0,
$$

$$
\det \frac{\partial^2 g}{\partial x^2} /(\vec{x} = 0) > 0. \eqno(7.4)
$$

 \ We explain in detail:

$$
 \frac{\partial g}{\partial x} = \left\{  \frac{ \partial g}{ \partial x_j} \right\} = \grad g,  \hspace{5mm}
\frac{\partial^2 g}{\partial x^2} = \left \{ \frac{\partial^2 g}{\partial x_k \partial x_l}  \right\}, \ -
$$
be a Hess matrix,  $ \ i,k,l = 1,2,\ldots,d. $ \par

 \ We assume  finally

$$
\lim_{|x| \to \infty} \frac{\partial^d g }{\prod_{k=1}^d  \partial x_k} = \infty.
$$

 \ We will denote the set of all such a functions by $  Y = Y(R^d) $ and denote also by $ D  $ introduced before matrix

$$
D = D_g := \frac{1}{2} \left\{ \frac{\partial^2 g(0)}{\partial x_k \partial x_l} \right\}.\eqno(7.5)
$$

 \ Evidently, the matrix $ D = D_g $ is non-negative definite,  write $  D = D_g \ge \ge 0.  $ \par

\vspace{3mm}

 \ {\bf  Definition 7.5 } \par

\vspace{3mm}

   \ Let the function  $  \phi = \phi(\lambda), \ \lambda \in R^d  $ be the Young function.
  We will say by definition likewise the one-dimensional case  that the {\it centered} (mean zero) random vector (r.v)
 $ \xi = \xi(\omega) = \vec{\xi} = (\xi(1), \xi(2), \ldots, \xi(d)) $ with values in the space $  R^d $ belongs to the
 space $ B(\phi), $  write $ \xi = \vec{\xi}\in  B(\phi),  $  if there exists certain non-negative constant $ \tau \ge 0 $ such that

$$
\forall \lambda \in R^d \ \Rightarrow
\max_{\vec{\epsilon}} {\bf E} \exp \left( \sum_{j=1}^d \epsilon(j) \lambda(j) \xi(j) \right) \le
\exp[ \phi(\lambda \cdot \tau) ]. \eqno(7.6)
$$

 \ The minimal value $ \tau $ satisfying (7.6) for all the values $  \lambda \in R^d, $
is named by definition as a $ B(\phi) \ $ norm of the vector $ \xi, $ write

$$
||\xi||B(\phi)  \stackrel{def}{=}
$$

 $$
 \inf \left\{ \tau, \ \tau > 0: \ \forall \lambda:  \ \lambda \in R^d \ \Rightarrow
 \max_{\vec{\epsilon}}{\bf E}\exp \left( \sum_{j=1}^d \epsilon(j) \lambda(j) \xi(j) \right) \le
 \exp(\phi(\lambda \cdot \tau)) \right\}. \eqno(7.7)
 $$

\vspace{3mm}

 \ The space $ \ B(\phi) \ $  relative introduced here norm $ ||\xi||B(\phi)  $ and ordinary algebraic operation
 is also multidimensional rearrangement invariant (symmetric) Banach space. \par

\vspace{3mm}

 \ For example, the {\it Moment Generating Function, briefly MGF, } $  \phi_{\xi}(\lambda) $ in these spaces, for the r.v. $ \ \xi \ $
 may be defined by the following natural way:

$$
\exp[ \phi_{\xi}(\lambda ) ] \stackrel{def}{=}
\max_{\vec{\epsilon} \in \Theta} {\bf E} \exp \left( \sum_{j=1}^d \epsilon(j) \lambda(j) \xi(j) \right),
 \eqno(7.8)
$$
if of course the random vector  $  \xi $ is centered and has an exponential tail of distribution. This imply  that
the natural function  $ \phi_{\xi}(\lambda ) $ is finite on some non-trivial central symmetrical neighborhood of origin,
or equivalently  the mean zero random vector  $  \xi $ satisfies the multivariate Kramer's condition. \par

 \ Obviously, for the natural function $ \phi_{\xi}(\lambda )  $

$$
||\xi||B(\phi_{\xi}) = 1.
$$

 \ It is easily to see that this choice of the generating function $ \phi_{\xi} $ is optimal, but in the practical using
often this function can not be calculated in explicit view, if of course there is a possibility to estimate its. \par

 \vspace{3mm}

 \ Note that the expression for the norm $ ||\xi||B(\phi)  $ dependent aside from the function $  \phi  $
only on the distribution $  \Law(\xi). $ Thus, this norm and correspondent space  $  B(\phi)  $ are rearrangement
invariant (symmetric) in the terminology of the classical book  \cite{Bennet1}, see chapters 1,2. \par

\vspace{4mm}

\ {\bf  I. The interrelations between MGF and tail behavior for the random vector. } \par

\vspace{4mm}

 \ The following important facts about upper tail estimate for the random vectors is proved in the preprint \cite{Ostrovsky6}. \par

\vspace{4mm}

\ {\bf  Corollary 7.1.a. } \ Let $  \phi = \phi(\lambda), \ \lambda \in R^d $ be arbitrary non-negative real valued function,
which is finite on some  non-empty symmetrical neighborhood of origin. Suppose for given centered $  d  \ - $ dimensional
random vector $ \xi = \vec{\xi}  $

$$
{\bf E} e^{(\lambda,\xi)} \le e^{ \phi( \lambda) }, \ \lambda \in R^d.  \eqno(7.9)
$$

 \ On the other words, $ ||\xi||B(\phi) \le 1. $ \ Then for all the non-negative deterministic vector $ x = \vec{x}  $
there holds

 $$
 U(\vec{\xi}, \vec{x}) \le \exp \left( - \phi^*(\vec{x}) \right) \ - \eqno(7.10)
 $$
the multidimensional generalization of Chernov's inequality. \par

\vspace{4mm}

 \ Moreover, the last estimate is essentially non-improvable, still in the one-dimensional case; there are some lower estimates
 for the considered here multivariate tail function in \cite{Ostrovsky6}. \par

\vspace{4mm}

 \ Let us obtain the converse conclusion. We retain the notations of lemma 3.1, namely $ \ K(\epsilon), \ \overline{K}(\lambda), $
 where $  X = R^d_+;  \  \zeta^{*'}(\lambda) = \grad \zeta^*(\lambda), $

$$
\theta(\lambda) \stackrel{def}{=} \frac{C_1}{ ( \lambda, \zeta^{*'}(\lambda) )} \eqno(7.11)
$$
for the greatest values $  |\lambda|:  \ |\lambda| > |\lambda_0|, $ where
$ \ \lambda_0 = \const \in R^d, \ |\theta(\lambda_0)| \le 0.5 $ (say); and apply the assertion of lemma 3.2.\par

\vspace{4mm}

\ {\bf Corollary 7.1.b. } Given, as above, for the {\it centered} r.v. $  \ \xi \ $

 $$
 U_{\xi}(x) \le \exp (  - \zeta(x)), \ x = \vec{x} \ge 0, \eqno(7.12)
 $$
where $  \zeta = \zeta(x) $ is suitable continuous non-negative function. \par

\ It is required to estimate for the sufficiently greatest values $ \ \lambda, $  say, $ |\lambda| > e,  $
 the moment generating function (MGF) for the r.v. $  \xi: $

$$
g_{\xi}(\lambda) \stackrel{def}{=} {\bf E} e^{ (\lambda, \xi) }, \ \lambda = \const \in R^d.
$$

\vspace{3mm}

 \ {\it   Both the theorems: theorem 3.2 and theorem 3.3 remains true, with alike proof,
  under formulated above multivariate notations and conditions. } \par

\vspace{3mm}

 \ The last statement may be reformulated as follows.  Assume as above the function $  \phi(\cdot)  $ be from the Young-Orlicz set,
 satisfying the restriction of corollary 7.2. The centered non-zero random vector $  \xi  $ belongs to the space  $   B(\phi): $

$$
\exists C_1 \in (0,\infty), \ \forall \lambda \in R^d
 \Rightarrow  {\bf E} e^{(\lambda,\xi)} \le e^{ \phi( C_1 \cdot \lambda) }, \ \lambda \in R^d,
$$
if and only if

$$
\exists C_2 \in (0,\infty), \ \forall \ x \in R^d_+  \Rightarrow \  U(\vec{\xi}, \vec{x}) \le \exp \left( - \phi^*(\vec{x}/C_2) \right).
$$

 \ More precisely, the following implication holds: there is a finite positive constant $  C_3 =  C_3(\phi)  $ such that
  for arbitrary non-zero centered r.v. $  \xi: \ ||\xi|| = ||\xi||B(\phi)  < \infty \ \Leftrightarrow $

$$
 \forall \lambda \in R^d
 \Rightarrow  {\bf E} e^{(\lambda,\xi)} \le e^{ \phi( ||\xi|| \cdot \lambda) }
$$
iff

$$
\exists C_3(\phi) \in (0,\infty) \ \forall \ x \in R^d_+  \Rightarrow \  U(\vec{\xi}, \vec{x}) \le
\exp \left( - \phi^*(\vec{x}/(C_3 / ||\xi||) \right).
$$

\vspace{4mm}

 \ {\bf Corollary 7.1.c.} Assume the non-zero centered random vector $  \xi = (\xi(1), \xi(2), \ldots, \xi(d)) $ belongs
to the space $  B(\phi):  $

$$
{\bf E} e^{(\lambda,\xi)} \le e^{ \phi(||\xi|| \cdot \lambda) }, \ \phi \in Y(R^d),
$$
and let $  y  $ be arbitrary positive non-random number. Then $  \forall y  > 0 \ \Rightarrow  $

$$
{\bf P} \left( \min_{j = 1,2,\ldots,n} |\xi(j)| > y  \right) \le 2^d \cdot
\exp \left(- \phi^*(y/||\xi||,y/||\xi||, \ldots, y/||\xi||) \right).  \eqno(7.13)
$$

 \vspace{3mm}

 \ The last estimate may be used in the analyse of discontinuous random fields, see for example \cite{Bickel1}, \cite{Billingsley1},
\cite{Billingsley2}, \cite{Gikhman1}, \cite{Grigorjeva1},  \cite{Prokhorov2}, \cite{Skorokhod1}. \par

\vspace{4mm}

{\bf Example 7.1.} Let as before $  V = R^d  $ and $  \phi(\lambda) = \phi^{(B)}(\lambda)  = 0.5(B\lambda,\lambda),  $ where
$  B  $  is non-degenerate positive definite symmetrical matrix, in particular $  \det B > 0. $ It follows from theorem 6.1
that the (centered) random vector $ \xi $ is subgaussian relative the matrix $ B: $

$$
\forall \lambda \in R^d \ \Rightarrow {\bf E} e^{ (\lambda, \xi)} \le e^{0.5 (B \lambda, \lambda) ||\xi||^2 }.
$$
iff for some finite positive constant $  K = K(B,d) $ and for any non-random positive vector $  x = \vec{x} $

$$
T_{\xi}(x) \le e^{- 0.5 \ \left( (B^{-1}x,x)/(K||\xi||^2) \right) }.
$$

\vspace{4mm}

\ {\bf  II. The interrelations between MGF and moments for the random vector. } \par

\vspace{4mm}

 \ This case is more complicated than considered before. We intend to generalize the results obtained in  \cite{Ostrovsky5},
 \cite{Ostrovsky6}, especially for the "inverse"  assertion. \par

 \ We need  getting to the presentation of the multidimensional case to extend our notations and restrictions.
 In what follows in this section the variables $ \lambda, r, x, \xi  $  are as before vectors from the space $  R^d, \ d = 2,3,\ldots, $
and besides $  r = \vec{r} = ( \ r(1), r(2), \ldots, r(d) \ ), \ r(j) \ge 1.  $\par

 \  A standard vector notations. Let $ a = \const \in R, $ then

$$
\vec{a}:= (a,a,\ldots, a), \ \dim \vec{a} = d.
$$

 \ Further,

 $$
  |r| = |\vec{r}| = \sum_j r(j), \hspace{4mm} |\xi| = |\vec{\xi}| =  \sqrt{ \ (\xi, \xi) \ };
 $$

$$
k! = \vec{k}! = \prod_{j=1}^d k(j)!, \ \vec{k} = (k(1), k(2), \ldots, k(d));
$$

$$
 \vec{x} \ge \vec{y} \ \Leftrightarrow \forall j \hspace{3mm} x(j) \ge  y(j);
$$

$$
 \vec{x} < \vec{y} \ \Leftrightarrow \forall j \hspace{3mm} x(j) <  y(j);
$$

 $$
 x^r = \vec{x}^{\vec{r}} = \prod_{j=1}^d x(j)^{r(j)}, \ \vec{x} \ge 0,
 $$

$$
\ln \vec{\lambda} = \{  \ln \lambda(1), \ \ln \lambda(2), \ldots, \ \ln \lambda(d) \}, \hspace{3mm} \vec{\lambda} > 0,
$$

$$
e^{\vec{\mu}} = \{ e^{\mu(1)}, \ e^{\mu(2)}, \ldots, \  e^{\mu(d)} \},
$$

$$
\Phi(\mu) = \Phi(\vec{\mu}) = \phi \left(e^{\vec{\mu}} \right),
$$

 $$
 \frac{r}{\lambda \cdot  e}  =  \frac{\vec{r}}{ \vec{\lambda} \ e} = \prod_{j=1}^d \left( \frac{r(j)}{ e \lambda(j)} \right) =
 e^{-|r|} \cdot \prod_{j=1}^d \left( \frac{r(j)}{\lambda(j)} \right),
 $$

$$
|\xi|_r = |\vec{\xi}|_{\vec{r}} = \left( {\bf E}|\vec{\xi}|^{\vec{r}} \right)^{1/|r|}.
$$

 \ We will use now the following elementary inequality

$$
x^r \le \left( \frac{r}{\lambda \ e} \right)^r \cdot e^{(\lambda, \ x)}, \hspace{3mm} r, \ \lambda, \ x > 0.
$$

 \ As a consequence: let $  \xi $ be non-zero  $ d  - $  dimensional mean zero random vector belonging to the space $  B(\phi).  $
Then

$$
{\bf E} |\xi|^r \le 2^d \  \left( \frac{r}{\lambda \ e} \right)^r  \ e^{ \phi(\lambda ||\xi||) } =
2^d \ e^{-|r|} \ r^r \ \lambda^{-r} \ e^{ \phi(\lambda ||\xi||) }, \ \lambda > 0.
$$

 \ We find likewise the one-dimensional case: \\

\vspace{4mm}

{\bf Proposition 7.2.a.} Let $ \phi (\cdot) $  be arbitrary non-negative continuous function and let the centered numerical
r.v. $  \xi $  be such that $  \xi \in B(\phi): \ 0 < ||\xi|| = ||\xi||B(\phi) < \infty. $ \par

 \ Then

$$
| \vec{\xi}|_{\vec{r}} \le e^{-1} \cdot  2^{d/|r|} \cdot \prod_j r(j)^{r(j)/|r|}  \cdot  e^{- \Phi^*(r)/|r| } \cdot ||\xi||B(\phi),
 \ r = \vec{r} > 0.
$$

 \ Note that in general case the expression $ |\xi|_r  $ does not represent the norm relative the random vector $  \vec{\xi}. $ \par
 \ But if we denote

$$
\psi_{\Phi}(\vec{r}) :=  e^{-1} \cdot 2^{d/|r|} \cdot \prod_j r(j)^{r(j)/|r|}  \cdot  e^{- \Phi^*(r)/|r| }
$$

  and define
$$
||\xi||G\psi_{\Phi} \stackrel{def}{=} \sup_{\vec{r} \ge 1} \left[  \frac{| \vec{\xi}|_{\vec{r}} }{\psi_{\Phi}(\vec{r}) }  \right], \eqno(7.14)
$$
we obtain some modification of the one-dimensional Grand Lebesgue Space (GLS) norm. \par

 \ The statement of proposition (7.2.a) may be rewritten as follows.

$$
||\xi||G\psi_{\Phi}  \le ||\xi||B(\phi). \eqno(7.15)
$$

 \vspace{3mm}

 \ Let us state the inverse up to multiplicative constant inequality. \par

\vspace{4mm}

 \ Given: for the mean zero random vector $ \xi = \vec{\xi} $

$$
|\vec{\xi}|_{\vec{r}} \le \psi^{|r|}( |r|) \cdot  e^{-\Delta (r)/|r|}, \eqno(7.16)
$$
or equally

\vspace{3mm}

$$
 {\bf E}|\vec{\xi}|^{\vec{r}} \le   |r|^{|r|} \cdot  e^{-\Delta (r)}. \eqno(7.16a)
$$

 \ We have for $  \vec{\lambda} \ge \vec{e} $

$$
{\bf E}e^{(\lambda, \xi)} - 1 \le  C(d) \ \Sigma \ldots \Sigma_{\vec{k} \ge \vec{2}} \
\left[\frac{ {\bf E} |\xi|^{\vec{k}}}{\vec{k}!} \right] \le
$$

$$
C_2(d) \ \Sigma \ldots \Sigma_{\vec{k} \ge \vec{2}} \left[ \ \vec{\lambda}^{\vec{k}} \cdot e^{- \Delta(\vec{k})} \ \right].
$$

 \ We have as above for the great values $ \lambda $ using again the Stirling's formula

$$
{\bf E} \exp (\lambda \xi) - 1 \le C_3(d) \ \sum_{\vec{k}\ge \vec{2}} \exp( \vec{k} \ln(C\lambda) - \Delta(\vec{k}))=
$$

$$
C_3(d) \ \sum_{\vec{k}\ge \vec{2}} \exp( \vec{k} \ln( \mu) - \Delta(\vec{k})), \hspace{5mm} \mu  = \vec{\mu}:= \ln (C_4 \vec{\lambda}).
$$

  \ Let $  \epsilon = \const \in (0,1); $ we apply the Young's inequality:

$$
k \mu \le \Delta(\epsilon k) + \Delta^*(\mu/\epsilon)
$$
and denoting as before

$$
\sigma_{\Delta}(\epsilon) = \sigma(\epsilon) := \sum_{\vec{k} \ge \vec{2}} e^{ \Delta(\epsilon k) - \Delta(k)}:
$$

$$
{\bf E}\exp(\lambda \xi) - 1 \le e^{\Delta^*(\mu/\epsilon)} \sum_{\vec{k} \ge \vec{2}} e^{ \Delta(\epsilon k) - \Delta(k)}  \le
\sigma_{\Delta}(\epsilon) \cdot e^{\Delta^*(\mu/\epsilon)}.\eqno(7.17)
$$

 \ Introduce also the function $  \Delta_1^*(\mu) $ as follows:

$$
 \Delta_1^*(\mu) := \inf_{\epsilon \in (0,1)} \ln \left\{ \sigma(\epsilon) \ \exp \left( \Delta^*(\mu/\epsilon)  \right) \right\};
 \eqno(7.18)
$$
then we obtain for the values $ \vec{\lambda} \ge \vec{e}: $

$$
{\bf E} e^{ (\lambda \xi)} \le e^{ \Delta_1^*(\ln (C_2 \lambda) )}.
$$

 \ At the same estimate is true also with another but again finite constant $ C_3 $ for the values  $ \vec{\lambda}  < - \vec{e}  $
and $ |\lambda| \le  e, $ following

$$
 \exists C_5 < \infty \ \Rightarrow  {\bf E} e^{(\lambda \xi)} \le e^{ \Delta_1^*(\ln (C_6 \ | \lambda)| )}.  \eqno(7.19)
$$

\vspace{3mm}

 \ Let us impose the following important condition $ (\Delta) $ on the source function $  \phi_{\Delta}(\cdot):  $

$$
 (\Delta): \hspace{5mm}   \exists C_7 \in [1, \infty) \Rightarrow \Delta_1^*(\ln (\ | \lambda)| ) \le
 \Delta^*(\ln (\ | C_8 \lambda)| ), \ \lambda = \vec{\lambda} \in R^d.  \eqno(7.20)
$$

\vspace{4mm}

 \ {\bf Proposition  7.2.b.} Suppose the function $ \psi = \psi_{\Delta}  $  satisfies the
condition $ \ (\Delta). \ $ Define the following function

$$
\phi_{\Delta}(\lambda) = \Delta^*(\ln |\lambda|), \ \lambda \in R.
$$

 \ Our proposition:

$$
|| \ \xi \ ||B \left(\phi_{\Delta} \right) \le C_4 \ || \ \xi \ ||G\psi_{\Delta}.  \eqno(7.21)
$$

\vspace{4mm}

{\bf III. \ Norm equivalence. }

\vspace{4mm}

 \ We define a function $  \psi(\cdot) $  in the form

$$
\psi_{(\phi)}(r) = \psi_{(\phi)}(\vec{r}) = |r| \cdot \exp \left[ - \beta^*(r)/|r|  \right], \ \vec{r} \ge \vec{1},
\eqno(7.22)
$$
where $  \beta(\cdot) $ is certain even  convex continuous function from the set $ \ \Phi. \   $

\vspace{3mm}

{\bf Proposition 7.2.c.} We have under formulated before conditions and notations

$$
|| \ \xi \ ||G\psi_{(\phi)} \le e^{-1} \ ||\ \xi \ ||B(\phi). \eqno(7.23)
$$

\vspace{4mm}

\ {\bf  III. The interrelations between tail function and moments norm for the random vectors. } \par

\vspace{4mm}

 \ {\bf A.  \ "Direct estimate".} \ Given: the random vector $  \  \xi = \vec{\xi} \ $ such that for some multivariate function
 $  \ \psi(\cdot) \in \Psi = \Psi_{\infty} \hspace{4mm} || \vec{\xi}|| =  ||\xi|| = ||\xi||G\psi \in (0,\infty). $  It is required
 to estimate the tail function $  \ U_{\xi}(x), $ for sufficiently greatest values $ x = \vec{x}, $ say $  \ \vec{x} > \vec{e}. \ $ \par

 \ Define the auxiliary function

$$
\nu(p) =  \nu_{\psi}(p) := |p| \cdot \ln \psi(p), \ p = \vec{p} \ge  \vec{1}
$$
and correspondingly

$$
\nu^*(z) = \nu^*_{\psi}(z) = \sup_{\vec{p} \ge \vec{1} } \left( (p,z) - \nu_{\psi}(p) \right), \ z \in R^d, \vec{z} \ge \vec{0}.
$$

 \vspace{3mm}

 \  Suppose $ \ \vec{\xi} \in G\psi, \ \xi \ne 0. \ $  We derive alike the second section, theorem 2.1, using again the Tchebychev -
 Markov inequality

$$
U_{\xi}(\vec{y}) \le \exp \left\{ \ - \nu^*[\ \ln (\vec{y}/||\xi||) \ ] \ \right\}, \ \vec{y} \ge \vec{e}. \eqno(7.23a)
$$

 \vspace{3mm}

 \ {\bf B. \ Let us deduce the "opposite" estimate.} \ Given:

 $$
 U_{\vec{\xi}}(\vec{x}) \le \exp (  - \zeta(\vec{x})), \ x \ge 0, \ \zeta(\vec{x}) \ge 0; \eqno(7.24)
 $$
and we denote $ Z(y) = \zeta(\exp y), \ y \in R^d.  $ \par

\vspace{3mm}

  \ Suppose $ Z(\cdot) $ is twice continuous differentiable
convex function on certain "octane"  $ \ Q_d(C) := (C, \infty)^d, \ C = \const > 0.  $
Denote by $ \Lambda = \Lambda(y) = \Lambda_Z(y)  $
the {\it minimal } eigenvalue of the (square) matrix $ Z^{''} (y) = \partial^2 Z/\partial y^2,   $
more detail

$$
  \left\{ \ Z^{''} (y) \ \right\}_{i,j} = \frac{\partial^2 Z}{\partial y_i \ \partial y_j}, \ = 1,2,\ldots d:
$$

$$
 \Lambda(y) = \Lambda_Z(y) := \inf_{x \in R^d, \ ||x||=1 } (Z^{''} (y) x,x),
$$

and suppose

$$
\exists C = \const > 0 \ \Rightarrow C_1:= \inf_{y \in Q_d(C)}  \Lambda_Z(y) > 0. \eqno(7.25)
$$
 Then

$$
|\xi|_p \le C_3 \exp \left(  Z^*(p)/p \right), \eqno(7.26)
$$
or equally

$$
||\xi||G(Z^*(p)/p) \le C_3 < \infty. \eqno(7.26.a)
$$

\vspace{3mm}

 \ Indeed, we have after integration "by parts"

$$
{\bf E} |\vec{\xi}|^{\vec{p}} \le  \left| \prod_j p_j \ \int_{R^d_+} \prod_j x_j^{p_j-1} \ U_{\xi}( \vec{x}) \ dx \right| \le
$$

$$
\le  \prod_j p_j \cdot \int_{R^d} e^{ (p,y) - Z(y)} dy =:  \prod_j p_j \ J(p).
$$

 \ We deduce:

$$
J(p) \le 2^d  \int_{ Q_d(C) } e^{(p,y) - Z(y)} \ dy.
$$

 \ Denote $  S = S(p,y) = (p,y) - Z(y) $ and $  y_0 = y_0(p) = \argmax_{y \ge C_1} S(p,y);  $ then

$$
S(p,y) \le \max_y S(p,y) - 0.5  \Lambda_Z(y) \ ||y - y_0||^2 \le
$$

$$
Z^*(p) - C_6 ||y - y_0||^2, \ C_6 = \const > 0,
$$
therefore

$$
J(p) \le C_7(d,Z)^p \ \exp \left( Z^*(p) \right), \ \vec{p} \ge \vec{1}. \eqno(7.27)
$$

 \ We used again the obvious estimate $  p^{1/p} \le e^{1/e}, \ p \ge 1. $ \par

 \vspace{4mm}

{\bf C. Coincidence.} \par

\vspace{4mm}

 \ It is convenient for us to rewrite the restriction (7.22) in the following form

$$
|\xi|_p \le |p| \ \ e^{ \ - \nu(p)/p \ }, \ p \ge 1, \eqno(7.28)
$$
i.e. in (7.14)

$$
\psi(p) := \psi_{\nu}(p) =  \ |p| \cdot \ e^{ - \nu(p)/p}.\eqno(7.28a)
$$

\vspace{4mm}

 \ We find analogously theorem 2.3 the following statement. \par

 \ Suppose that the function $  \nu = \nu(p), \ p = \vec{p} \ge \vec{1}  $ in (7.28) is continuous, convex, and
such that the function $  y \to \nu(\exp y) $ satisfies the condition (6.25). Then the GLS  vector norm estimate for the
non - zero r.v. $ \xi $ of the form

$$
|\xi|_p \le C_1 \ |p| \cdot e^{ \ - \nu(p)/p \ }, \ p \ge 1, \eqno(7.29)
$$
is quite equivalent to the following tail inequality

$$
U_{\vec{\xi}}(\vec{y}) \le \exp \left( - \nu^*(\ln (\vec{y}/C_2)  \right), \
C_2 = \const \in (0,\infty), \ \vec{y} \ge C_2 \vec{e}. \eqno(7.30)
$$

\vspace{4mm}

\section{Multidimensional example}

\vspace{4mm}

 \ We can construct many multivariate examples by the following way. Recall that the function $ \ g = g(x), \ g: R^d \to R $
 is said to be {\it radial},  or {\it spherical symmetric, } if it dependent only on the Euclidean length of the variable $  x: $

$$
\exists g_1: R_+ \to R, \ g(x) = g_1(|x|)  = g_1 ( (x,x)^{1/2} ). \eqno(7.1)
$$

 \ It is easily to verify that the Young-Fenchel transform of radial function is also radial function. Indeed, let the function
$ g = g(x), \   x \in R^d  $  be spherical symmetric, i.e. satisfies (7.1.) Let also $ \ A \ $ be arbitrary orthogonal matrix:
$ \ A \ A^T = E. \  $ We deduce

$$
g^*(y) = \sup_{x \in R^d} ( (x,y) - g_1(|x|) ) = \sup_{z \in R^d}((Az,y) - g_1(|z|) )  =
$$

$$
\sup_z ( (z, A^Ty) - g_1(|z|) )  = g_1^*(|A^T y|) = g_1^*(| y|).
$$

\vspace{3mm}

 \ As a slight consequence: for any random  vector $ \eta = \vec{\eta} $ the estimate of the form

$$
{\bf E} \exp(\lambda,\eta) \le \exp(C_1 |\lambda|^m), \ |\lambda| \ge 1, \ m = \const > 1 \eqno(7.2)
$$
is quite equivalent to the tail estimate

$$
U_{\vec{\eta}}(x) \le \exp \left( - C_2 |x|^{m/(m-1)} \right), \ \min_j |x_j| \ge 1.\eqno(7.3)
$$

 \vspace{3mm}

 \ Note that if the function $ \phi = \phi(\mu ), \ \mu \in R  $ is from the set $ \Phi,  $ then the
 multivariate function

$$
\chi (\vec{\lambda}) = \phi(|\vec{\lambda} | )
$$
 can serve as an example of the multivariate function belonging  to the  set $ \Phi: R^d \to R. $ \par

 \ Analogously may be considered the case of a functions of the form

$$
\chi (\vec{\lambda}) = \phi( (A \ \vec{\lambda}, \vec{\lambda})^{1/2} ),
$$
where $  A  $ is symmetrical positive definite matrix: $  A \in R^d \otimes  R^d. $\par

\vspace{4mm}

\section{Concluding remarks}

\vspace{4mm}

 \hspace{4mm} {\bf 1.} Note that the obtained in this report results does not follow from ones in the articles
 \cite{Ostrovsky201}, \cite{Ostrovsky202}, \cite{Ostrovsky203}. Considered here random variables have the exponential
decreasing tails of distributions, in contradiction to the considered in mentioned reports. \par

\vspace{3mm}

 \ {\bf 2.} Open problem: under which additional conditions imposed on the function $ \psi = \psi(p)  $
  the assertion of theorem 3.3 remains true in the case when $  \ b < \infty ? $ \par

\vspace{3mm}

 \ {\bf 3.} It looks like a multi-dimensional case when the support of the function $  \psi(\cdot) $ is (partially) bounded? \par

\end{document}